\documentclass[11pt,a4paper]{article}
\usepackage{amsfonts}
\usepackage{amsmath}
\usepackage{amssymb}
\usepackage{amsthm}
\usepackage{color}

\newtheorem{theorem}{Theorem}
\newtheorem{thma}{Theorem}

\newtheorem{corollary}[theorem]{Corollary}

\newtheorem{lemma}[theorem]{Lemma}

\newtheorem{proposition}[theorem]{Proposition}

\begin{document}

\title{Analytic content and the isoperimetric inequality in higher dimensions%
}
\date{}
\author{Stephen J. Gardiner, Marius Ghergu\ and Tomas Sj\"{o}din}
\maketitle

\begin{abstract}
This paper establishes a conjecture of Gustafsson and Khavinson, which
relates the analytic content of a smoothly bounded domain in $\mathbb{R}^{N}$
to the classical isoperimetric inequality. The proof is based on a novel
combination of partial balayage with optimal transport theory.
\end{abstract}

\section{Introduction}

\footnotetext{%
\noindent 2010 \textit{Mathematics Subject Classification } 31B05.
\par
\noindent \textit{Keywords: }analytic content, harmonic vector fields,
isoperimetric inequality, optimal transport, partial balayage
\par
{}
\par
\noindent}Let $\omega $ be a bounded domain in the complex plane $\mathbb{C} 
$ such that $\partial \omega $ is the disjoint union of finitely many simple
analytic curves, and let $\mathcal{A}(\omega )$ denote the collection of
continuous functions on $\overline{\omega }$ that are analytic on $\omega $.
Further, let $\left\Vert g\right\Vert _{S}$ denote $\sup_{S}\left\vert
g\right\vert $ for any bounded function $g:S\rightarrow \mathbb{C}$. The 
\textit{analytic content} of $\omega $ is then defined by 
\begin{equation*}
\lambda (\omega )=\inf \{\left\Vert \overline{z}-\phi \right\Vert _{%
\overline{\omega }}:\phi \in \mathcal{A}(\omega )\}.
\end{equation*}%
The inequalities for $\lambda (\omega )$ given below, which imply the
classical isoperimetric inequality, are due to Alexander \cite{Al} and
Khavinson \cite{Kh84}.

\begin{thma}
Let $A$ and $P$ denote the area and perimeter of $\omega $, respectively.
Then 
\begin{equation*}
\frac{2A}{P}\leq \lambda (\omega )\leq \sqrt{\frac{A}{\pi }}.
\end{equation*}
\end{thma}

\bigskip

An exposition of this circle of ideas may be found in Gamelin and Khavinson 
\cite{GaKh}, and a wider survey of related results is provided by B\'{e}n%
\'{e}teau and Khavinson \cite{BK}. It was shown in \cite{GaKh} that equality
with the upper bound occurs if and only if $\omega $ is a disc. Recently,
Abanov et al \cite{ABKT} have shown that equality with the lower bound
occurs if and only if $\omega $ is a disc or an annulus.

Rewriting $\lambda (\omega )$ as $\inf \{\left\Vert z-\overline{\phi }%
\right\Vert _{\overline{\omega }}:\phi \in \mathcal{A}(\omega )\}$, it can
be seen that a natural generalization of this quantity to smoothly bounded
domains $\Omega $ in Euclidean space $\mathbb{R}^{N}$ ($N\geq 2$) is given by%
\begin{equation*}
\lambda (\Omega )=\inf \{\left\Vert x-f\right\Vert _{\overline{\Omega }%
}:f\in A(\Omega )\},
\end{equation*}%
where $A(\Omega )$ denotes the space of \textit{harmonic vector fields }$%
f=(f_{1},...,f_{N})\in C(\overline{\Omega })\cap C^{1}(\Omega )$ and 
\begin{equation*}
\left\Vert f\right\Vert _{S}=\sup_{S}\left\Vert f\right\Vert ,\text{ \ where
\ }\left\Vert f\right\Vert =\sqrt{f_{1}^{2}+...+f_{N}^{2}}.
\end{equation*}%
(Thus $f$ satisfies $\mathrm{div}f=0$ and $\mathrm{curl}f=0$, where the
latter condition means that 
\begin{equation*}
\frac{\partial f_{j}}{\partial x_{k}}-\frac{\partial f_{k}}{\partial x_{j}}=0%
\text{ \ \ for all }j,k\in \{1,...,N\}\text{\ on\ }\Omega .)
\end{equation*}%
If we write 
\begin{equation*}
\mathcal{H}=\{h\in C^{1}(\overline{\Omega }):\Delta h=0\text{ \ on \ }\Omega
\},
\end{equation*}%
then 
\begin{equation}
\{\nabla h:h\in \mathcal{H}\}\subset A(\Omega ).  \label{B}
\end{equation}%
Let $r_{\Omega }>0$ be chosen so that a ball of radius $r_{\Omega }$ has the
same volume as $\Omega $. Gustafsson and Khavinson \cite{GK1} established
the following inequalities for $\lambda (\Omega )$ in higher dimensions.

\begin{thma}
\label{GK}Let $\Omega $ be a bounded domain in $\mathbb{R}^{N}$ ($N\geq 3$)
with volume $V$ such that $\partial \Omega $ is the disjoint union of
finitely many smooth components with total surface area $P$. Then there
exists a constant $c_{N}>1$ such that 
\begin{equation}
\frac{NV}{P}\leq \lambda (\Omega )\leq c_{N}r_{\Omega }.  \label{GKi}
\end{equation}
\end{thma}

\bigskip

The lower bound in (\ref{GKi}) is sharp, since $\lambda (\Omega )=r$ when $%
\Omega $ is a ball of radius $r$ (see Theorem 3.1 in \cite{GK1}). Regarding
the upper bound, Gustafsson and Khavinson conjectured that the constant $%
c_{N}$ may be replaced by $1$, in which case (\ref{GKi}) would again contain
the classical isoperimetric inequality. However, the methods of \cite{GK1}
do not yield such a conclusion. The purpose of this paper is to verify this
long-standing conjecture.

Following \cite{GK1} we define a related domain constant,%
\begin{eqnarray}
\lambda _{1}(\Omega ) &=&\inf \{\left\Vert x-\nabla h\right\Vert _{\overline{%
\Omega }}:h\in \mathcal{H}\}  \notag \\
&=&\inf \{\left\Vert \nabla u\right\Vert _{\overline{\Omega }}:u\in C^{1}(%
\overline{\Omega })\text{ \ and \ }\Delta u=N\text{ on }\Omega \},
\label{L1}
\end{eqnarray}%
where we have used the observation that a function $u$ in $C^{1}(\overline{%
\Omega })$\ satisfies\ $\Delta u=N$ on $\Omega $ if and only if the function 
$h$ defined by $h(x)=\left\Vert x\right\Vert ^{2}/2-u(x)$ belongs to $%
\mathcal{H}$. {For future reference we note that, for such $u$, it follows
from the harmonicity of the partial derivatives of $u$ that $\Vert \nabla
u\Vert $ is subharmonic on }$\Omega ${\ (cf. Theorem 3.4.5 of \cite{AG}),
and so 
\begin{equation}
\Vert \nabla u\Vert _{\overline{\Omega }}=\Vert \nabla u\Vert _{\partial {%
\Omega }}  \label{gradsubh}
\end{equation}%
by the maximum principle.}

It follows from (\ref{B}) that 
\begin{equation*}
\lambda (\Omega )\leq \lambda _{1}(\Omega ).
\end{equation*}%
(Equality holds when $\Omega $ is simply connected.) Gustafsson and
Khavinson actually showed in \cite{GK1} that $\lambda _{1}(\Omega )\leq
c_{N}r_{\Omega }$ for an explicit constant $c_{N}>1$. We will establish the
following estimate.

\begin{theorem}
\label{main}Let $\Omega $ be a bounded domain in $\mathbb{R}^{N}$ such that $%
\partial \Omega $ is the disjoint union of finitely many smooth components.
Then $\lambda _{1}(\Omega )\leq r_{\Omega }$. Further, equality holds if and
only if $\Omega $ is a ball.
\end{theorem}

In the light of Theorem \ref{GK} we immediately arrive at the following
conclusion.

\begin{corollary}
Let $\Omega $ be a bounded domain in $\mathbb{R}^{N}$ with volume $V$ such
that $\partial \Omega $ is the disjoint union of finitely many smooth
components with total surface area $P$. Then 
\begin{equation*}
\frac{NV}{P}\leq \lambda (\Omega )\leq r_{\Omega }.
\end{equation*}%
Further, $\lambda (\Omega )=r_{\Omega }$ if and only if $\Omega $ is a ball.
\end{corollary}

The proof of Theorem \ref{main} combines the technique of partial balayage
with results from the theory of optimal transport. Later we will discuss
separate necessary and sufficient conditions for a function $u$ in $C^{1}(%
\overline{\Omega })$ to be a minimizer for $\lambda _{1}(\Omega )$ in (\ref%
{L1}), whenever such minimizers exist.

The purpose of Corollary 2 is to relate analytic content to the
isoperimetric inequality. We do not claim that it offers a novel or shorter
proof of the latter, since (apart from other more classical proofs) there
are already proofs using optimal transport theory as in McCann and Guillen 
\cite{MCG}, and Cabr\'e{} \cite{Cab} had previously provided a short proof
of it based on more geometric methods.

\section{Proof of Theorem \protect\ref{main}}

\subsection{Tools for the proof\label{Tools}}

Let $m$ denote Lebesgue measure on $\mathbb{R}^{N}$ and $O\subset \mathbb{R}%
^{N}$ be a bounded domain. Further, let $G_{O}(\cdot ,\cdot )$ denote the
Green function of $O$, and $G_{O}\mu ,G_{O}\nu $ denote the potentials of
(positive) measures $\mu ,\nu $ on $O$, where $\nu \ll m$. The Green
function is normalized so that $-\Delta G_{O}\gamma =\gamma $ in the sense
of distributions for any potential $G_{O}\gamma $. We define 
\begin{equation}
P_{\mu }^{\nu }=G_{O}\nu +\sup \{s:s\text{ is subharmonic on }O\text{ and }%
s\leq G_{O}\mu -G_{O}\nu \}\text{ \ on }O,  \label{pb}
\end{equation}%
whence $P_{\mu }^{\nu }\leq G_{O}\mu $, and recall the following facts (see 
\cite{GSa} or \cite{GS}).

\begin{thma}
\label{PB}(a) $P_{\mu }^{\nu }=G_{O}\eta $ for some measure $\eta $ on $O$
satisfying $\eta \leq \nu $.\newline
(b) $\eta =\nu |_{S}+\mu |_{O\backslash S}$, where $S=\{G_{O}\eta <G_{O}\mu
\}$.
\end{thma}

We will refer to the measure $\eta $ in the above theorem as the (\textit{%
partial}) \textit{balayage of }$\mu $\textit{\ onto }$\nu $\textit{\ in }$O$%
, and denote it by $\mathcal{B}\mu $, where $\nu $ is to be understood from
the context. We note that $G_{O}\mu -G_{O}\eta $ is the smallest nonnegative
lower semicontinuous function $w$ on $O$ satisfying $-\Delta w\geq \mu -\nu $
in the sense of distributions. Thus, if $\mu _{1}\geq \mu $, then the set $%
S(\mu )$ associated with $\mu $\ is contained in the corresponding set $%
S(\mu _{1})$. It follows that 
\begin{equation}
\mu _{1}\geq \mu \Longrightarrow \mathcal{B}\mu _{1}\geq \mathcal{B}\mu .
\label{ineqB}
\end{equation}%
We will also need the following lemma. Let $B(x,r)$ denote the open ball in $%
\mathbb{R}^{N}$ with centre $x$ and radius $r$.

\begin{lemma}
\label{Lsig}Let $\nu =Nm|_{O}$ and $\overline{\Omega }\subset \Omega
_{0}\subset O$,\ where $\Omega _{0}$ is another open set. If $\tau $ is a
measure with $\mathrm{supp}\tau \subset \overline{\Omega }$, then there
exists $b>0$ such that 
\begin{equation*}
\mathrm{supp}\mathcal{B}(Nm|_{\Omega }+b\tau )\subset \Omega _{0}.
\end{equation*}
\end{lemma}

\noindent \textbf{Proof. }Let $\Omega ^{\prime }$ be an open set such that $%
\overline{\Omega }\subset \Omega ^{\prime }$ and $\overline{\Omega ^{\prime }%
}\subset \Omega _{0}$, and let 
\begin{equation*}
\varepsilon =2^{-1}\mathrm{dist}(\partial \Omega ^{\prime },\Omega \cup (%
\mathbb{R}^{N}\backslash \Omega _{0})).
\end{equation*}%
Let $\tau ^{\ast }$ be the sweeping (classical balayage) of $\tau $ onto $%
\partial \Omega ^{\prime }$, and define 
\begin{equation*}
\gamma (x)=\int_{\partial \Omega ^{\prime }}\phi _{\varepsilon }(x-y)d\tau
^{\ast }(y)\text{ \ \ \ }(x\in \mathbb{R}^{N}),
\end{equation*}%
where $\phi _{\varepsilon }$ is a non-negative rotationally invariant $%
C^{\infty }$ smoothing kernel on $\mathbb{R}^{N}$\ with support $\overline{B}%
(0,\varepsilon )$ (see, for example, Section 3.3 of \cite{AG}). We choose $b$
sufficiently small that $b\gamma \leq N$, whence 
\begin{equation*}
P_{G_{O}(Nm|_{\Omega }+b\gamma m)}^{\nu }=G_{O}(Nm|_{\Omega }+b\gamma m).
\end{equation*}%
Since $G_{O}(\gamma m)\leq G_{O}\tau ^{\ast }\leq G_{O}\tau $, with equality
outside $\{x\in O:\mathrm{dist}(x,\Omega ^{\prime })\leq \varepsilon \}$, we
see that 
\begin{equation*}
G_{O}(Nm|_{\Omega }+b\gamma m)\leq P_{Nm|_{\Omega }+b\tau }^{\nu }\leq
G_{O}(Nm|_{\Omega }+b\tau ),
\end{equation*}%
again with equality outside $\{x\in O:\mathrm{dist}(x,\Omega ^{\prime })\leq
\varepsilon \}$, and so 
\begin{equation*}
\mathrm{supp}\mathcal{B}(Nm|_{\Omega }+b\tau )\subset \{x\in O:\mathrm{dist}%
(x,\Omega ^{\prime })\leq \varepsilon \}.
\end{equation*}

\bigskip

We recall the following composite result from the theory of optimal
transport, in which the existence and smoothness of the function $v$ are due
to Brenier \cite{Br} and Caffarelli \cite{Ca}, respectively. (See also
Chapters 3 and 4 of Villani's book \cite{Vi}.)

\begin{thma}
\label{brenier} Let $D\subset \mathbb{R}^{N}$ be a bounded open set such
that $m(\partial D)=0$. Then there exists a convex function $v:\mathbb{R}%
^{N}\rightarrow (-\infty ,\infty ]$ which is $C^{2}$ on $D$, and for which $%
\nabla v$ maps $D$ into  $B(0,r_{D})$ and is
measure-preserving, in the sense that $m(A)=m((\nabla v)(A))$ for any Borel
set $A\subset D$.
\end{thma}

\begin{lemma}
\label{Lem}If a measure $\gamma $ on $\Omega $ has bounded density with
respect to $m$, then $G_{\Omega }\gamma \in C^{1}(\overline{\Omega })$, and 
\begin{equation}
\frac{\partial G_{\Omega }\gamma }{\partial y_{i}}(y)=\int_{\Omega }\frac{%
\partial G_{\Omega }}{\partial y_{i}}(x,y)d\gamma (x)\text{ \ \ \ }(y\in 
\overline{\Omega };i=1,...,N).  \label{under}
\end{equation}
\end{lemma}

To see this, we note that standard arguments (cf. Theorem 4.5.3 of \cite{AG}%
) show that $G_{\Omega }\gamma \in C^{1}(\Omega )$ and that (\ref{under})
holds when $y\in \Omega $. We now fix $i$ and note (see Widman \cite{Wi})
that $\left( \partial G_{\Omega }/\partial y_{i}\right) (x,\cdot )$\ has a
continuous extension to $\overline{\Omega }$ for each $x\in \Omega $.\ Let $%
y_{0}\in \partial \Omega $ and $\varepsilon >0$, and define 
\begin{equation*}
\psi _{j}(y)=\int_{\Omega _{j}}\frac{\partial G_{\Omega }}{\partial y_{i}}%
(x,y)d\gamma (x)\text{ \ \ \ }(y\in \overline{\Omega };j=1,2),
\end{equation*}%
where $\Omega _{1}=\Omega \backslash B(y_{0},\varepsilon )$ and $\Omega
_{2}=\Omega \cap B(y_{0},\varepsilon )$. Then $\psi _{1}$ is continuous at $%
y_{0}$, and (by estimates in \cite{Wi}) 
\begin{equation*}
\left\vert \psi _{2}(y)\right\vert \leq C(\Omega )\varepsilon \left\Vert 
\frac{d\gamma }{dm}\right\Vert _{L^{\infty }(\Omega )}\text{ \ \ \ }(y\in
B(y_{0},\varepsilon /2)\cap \overline{\Omega }).
\end{equation*}%
It follows that 
\begin{equation*}
\frac{\partial G_{\Omega }\gamma }{\partial y_{i}}(y)\rightarrow
\int_{\Omega }\frac{\partial G_{\Omega }}{\partial y_{i}}(x,y_{0})d\gamma (x)%
\text{ \ \ \ }(y\rightarrow y_{0}).
\end{equation*}

\subsection{Proof of the inequality\label{ineq}}

Let 
\begin{equation*}
b_{N}=\frac{m(\{y\in B(0,1):y_{N}\geq 1/2\})}{m(B(0,1))},
\end{equation*}%
and let $D$ be a bounded open set such that $\overline{\Omega }\subset D$, $%
m(\partial D)=0$ and $m(D\backslash \Omega )<b_{N}m(\Omega )$. We next
choose $v$ as in Theorem \ref{brenier}. Since $\nabla v$ is
measure-preserving on $D$, the Hessian of $v$, which is positive
semi-definite because $v$ is convex, has determinant equal to $1$, and so $%
\Delta v\geq N$ by the arithmetic-geometric means inequality for the
eigenvalues of the Hessian (cf. the argument in Section 1.6 of McCann and
Guillen \cite{MCG}). It will be enough to show that $\lambda _{1}(\Omega
)\leq r_{D}$, since $r_{D}$ can be made arbitrarily close to $r_{\Omega }$.
This inequality trivially holds if $\Delta v\equiv N$ on $\Omega $, so we
assume from now on that $(\Delta v-N)m|_{\Omega }\neq 0$.

Let $R$ be an open set satisfying $\overline{\Omega }\subset R$ and $%
\overline{R}\subset D$. Since 
\begin{equation*}
m(D\backslash R)\leq m(D\backslash \Omega )<b_{N}m(\Omega )<b_{N}m(D)
\end{equation*}%
and $\nabla v$ is measure preserving on $D$, we see that%
\begin{equation}
(\nabla v)(R)\cap \{y\in B(r_{D}):y\cdot x\geq r_{D}/2\}\neq \emptyset \text{
\ \ \ }(x\in \partial B(0,1)).  \label{half}
\end{equation}%
Also, since $v$ is convex, the function 
\begin{equation*}
w(x)=\sup \left\{ v(y)+\nabla v(y)\cdot (x-y):y\in R\right\} \text{ \ \ \ }%
(x\in \mathbb{R}^{N})
\end{equation*}%
equals $v$ on $R$. Clearly $w$ is subharmonic (and indeed convex) on $%
\mathbb{R}^{N}$.

We now define $\varepsilon =2^{-1}\mathrm{dist}(\overline{\Omega },\mathbb{R}%
^{N}\backslash R)$ and%
\begin{equation*}
w_{\varepsilon }(x)=\int_{B(\varepsilon )}\phi _{\varepsilon }(x-y)w(y)dm(y)%
\text{ \ \ \ }(x\in \mathbb{R}^{N}),
\end{equation*}%
where $\phi _{\varepsilon }$ is a smoothing kernel, as before. Then $%
w_{\varepsilon }$ is also subharmonic (and convex) on $\mathbb{R}^{N}$, and $%
w_{\varepsilon }\geq w$ (see Theorem 3.3.3 in \cite{AG}). We further define%
\begin{equation*}
O=\{x\in \mathbb{R}^{N}:w_{\varepsilon }(x)<C\},
\end{equation*}%
where%
\begin{equation*}
C>\sup \{w_{\varepsilon }(x):x\in R\}.
\end{equation*}%
The set $O$ clearly contains $\overline{R}$. It is also bounded, since for
any $x\in \partial B(0,1)$ we see from (\ref{half}) that there exists $%
y_{x}\in R$ such that $\nabla v(y_{x})\cdot x\geq r_{D}/2$, and so 
\begin{eqnarray*}
w_{\varepsilon }(tx) &\geq &w(tx)\geq v(y_{x})+\nabla v(y_{x})\cdot
(tx-y_{x}) \\
&=&t\nabla v(y_{x})\cdot x+v(y_{x})-\nabla v(y_{x})\cdot y_{x} \\
&\geq &r_{D}t/2+\inf \{v(y)-\nabla v(y)\cdot y:y\in R\}\text{ \ \ \ }(t>0).
\end{eqnarray*}%
By Sard's theorem and the smoothness of $w_{\varepsilon }$, we can choose $C$
such that the set $O$ is smoothly bounded and 
\begin{equation*}
\Omega _{0}\subset O\text{, \ where }\Omega _{0}=\{x\in \mathbb{R}^{N}:%
\mathrm{dist}(x,\Omega )<1\}.
\end{equation*}%
Since $w_{\varepsilon }=C$ on $\partial O$, the function $C-w_{\varepsilon }$
is the potential $G_{O}\mu $ of the measure $\mu =(\Delta w_{\varepsilon })m$
on $O$. Further, $\mu \geq Nm|_{\Omega }$, since $\Delta w=\Delta v\geq N$
on $R$, and so 
\begin{equation*}
\Delta w_{\varepsilon }(x)=\int_{B(\varepsilon )}\phi _{\varepsilon
}(x-y)(\Delta w)(y)dm(y)\geq N\int_{B(\varepsilon )}\phi _{\varepsilon
}(x-y)dm(y)=N\text{ \ \ \ }(x\in \Omega ).
\end{equation*}

We next apply Theorem \ref{PB} with $\nu =Nm|_{O}$. The partial balayage $%
\eta =\mathcal{B}\mu $ satisfies $\eta =Nm|_{S}+\mu |_{O\setminus S},$ where 
$S=\{G_{O}\eta <G_{O}\mu \}$. Since $G_{O}\mu \geq G_{O}\eta $ and $\mu
=(\Delta w_{\varepsilon })m\geq Nm=\eta $ on $\Omega $, the function $%
G_{O}\mu -G_{O}\eta $ is nonnegative and superharmonic on $\Omega $. In
fact, since $\Omega $ is connected and 
\begin{equation*}
(\mu -\eta )(\Omega )=\int_{\Omega }(\Delta v-N)dm>0,
\end{equation*}%
it is strictly positive there by the {minimum} principle, and so $\Omega
\subset S$. We note that $G_{O}\mu ,G_{O}\eta \in C^{1}(\overline{O})$, by
Lemma \ref{Lem}. Since the nonnegative function $G_{O}\mu -G_{O}\eta $
achieves its minimum value $0$ on $O\backslash S$, we have $\nabla G_{O}\mu
=\nabla G_{O}\eta $ on $\partial S\cap O$. Also, since $G_{O}\eta \leq
G_{O}\mu $, we have%
\begin{equation*}
\Vert \nabla G_{O}\eta \Vert =-\frac{\partial }{\partial n}G_{O}\eta \leq -%
\frac{\partial }{\partial n}G_{O}\mu =\Vert \nabla G_{O}\mu \Vert \text{ \ \
\ on \ }\partial O,
\end{equation*}%
where $n$ denotes the outward unit normal to $\partial O$. Thus%
\begin{equation*}
\Vert \nabla G_{O}\eta \Vert _{\partial S}\leq \sup \left\{ \Vert \nabla
G_{O}\mu (x)\Vert :x\in O\right\} =\sup \left\{ \Vert \nabla w_{\varepsilon
}(x)\Vert :x\in O\right\} \leq r_{D},
\end{equation*}%
because $w$ (and hence also $w_{\varepsilon }$) is Lipschitz on $\mathbb{R}%
^{N}$ with Lipschitz constant at most $r_{D}$. {Finally, since }$\Omega
\subset S$ and {$\eta =Nm$ in $S$, it follows (see (\ref{gradsubh})) that 
\begin{equation}
\lambda _{1}(\Omega )\leq \Vert \nabla G_{O}\eta \Vert _{\overline{\Omega }%
}\leq \Vert \nabla G_{O}\eta \Vert _{\overline{S}}=\Vert \nabla G_{O}\eta
\Vert _{\partial S}\leq r_{D},  \label{bound}
\end{equation}%
}as desired.

\subsection{The case of equality}

The following result strengthens the conclusion of the previous section.

\begin{proposition}
Let $v$ be as in Theorem \ref{brenier}, with $D=\Omega $. If $(\Delta
v-N)m|_{\Omega }\neq 0$, then there is a domain $U$ containing $\overline{%
\Omega }$ and a function $u\in C^{1}(U)$ satisfying $\Delta u=N$ and $%
\left\Vert \nabla u\right\Vert \leq r_{\Omega }$ in $U$.
\end{proposition}

\noindent \textbf{Proof. }We may assume that $v(x_{0})=0$ for some{\ $%
x_{0}\in \Omega $. Let 
\begin{equation*}
D_{k}=\{x\in \mathbb{R}^{N}:\mathrm{dist}(x,\Omega )<\delta _{k}\}\text{ \ \
\ }(k\geq 0),
\end{equation*}%
where $(\delta _{k})_{k\geq 0}$ is a strictly decreasing sequence of
positive numbers with limit $0$ and $\delta _{0}<1$ is chosen small enough
so that $m(D_{0}\backslash \Omega )<b_{N}m(\Omega )$. For each }$k$ we
choose $v_{k}$ as in Theorem \ref{brenier}, with $D=D_{k}$, and such that 
$v_{k}(x_{0})=0$. {W}e next choose open sets $R_{k}$ such that $\overline{%
\Omega }\subset R_{k}$ and $\overline{R_{k}}\subset D_{k}$, and define $%
\varepsilon _{k}=2^{-1}\mathrm{dist}(\overline{\Omega },\mathbb{R}%
^{N}\setminus R_{k})$.

For each $k$ we apply the construction of \S \ref{ineq} with $%
D=D_{k},R=R_{k},v=v_{k},$ and $\varepsilon =\varepsilon _{k}$. Propositions
3.1 and 3.2 of Brenier \cite{Br2}, applied to the measures 
\begin{equation*}
\frac{m\left\vert _{D_{k}}\right. }{m(D_{k})}\text{ \ }(k\in \mathbb{N})%
\text{ \ \ \ and \ \ \ }\frac{m\left\vert _{B(0,r_{\Omega })}\right. }{%
m(B(0,r_{\Omega }))},
\end{equation*}%
show that $v_{k}\rightarrow v$ uniformly on $\Omega $. 

The functions $w$ and $w_{\varepsilon }$ in \S \ref{ineq} will now be
denoted, by abuse of notation, $w_{k}$ and $w_{k,\varepsilon _{k}}$
respectively. {Since  $\Vert \nabla w_{k}\Vert \leq r_{D_{k}}\leq r_{D_{0}}$
on $\mathbb{R}^{N}$ by construction, we see that $|w_{k}-w_{k,\varepsilon
_{k}}|\leq r_{D_{0}}\varepsilon _{k}$ on $\mathbb{R}^{N}$. Hence }$\left( {%
w_{k,\varepsilon _{k}}}\right) ${\ converges uniformly to $v$ on $\Omega $,
in view of the fact that $v_{k}=w_{k}$ on $R_{k}$, which contains $\Omega $.}

We now choose numbers $C_{k}$ so that the {set} 
\begin{equation*}
O_{k}=\{x\in \mathbb{R}^{N}:w_{k,\varepsilon _{k}}(x)<C_{k}\}
\end{equation*}%
satisfies $D_{0}\subset O_{k}$ for each $k$. Since the sequences $\left(
\left\Vert v_{k}\right\Vert _{L^{\infty }(D_{k})}\right) $ and \newline $\left(
\left\Vert \nabla v_{k}\right\Vert _{L^{\infty }(D_{k})}\right) $ are
bounded, we can furthermore arrange that the set \newline $O=\cup _{k}O_{k}$ is
bounded.

Let $\gamma_{k}=(\Delta w_{k,\varepsilon_k})m|_{\Omega }$ and $\gamma
=(\Delta v)m|_{\Omega }$. These are non-negative measures, and the
divergence theorem shows that%
\begin{equation}
\Vert \gamma _{k}\Vert =\int_{\Omega }\Delta
w_{k,\varepsilon_k}dm=\int_{\partial \Omega }\frac{\partial
w_{k,\varepsilon_k}}{\partial n}d\sigma \leq r_{D_{k}}\sigma (\partial
\Omega )\text{ \ \ \ }(k\in \mathbb{N}),  \label{gamk}
\end{equation}%
where $\sigma $ denotes surface area measure on $\partial \Omega $. Since $%
\int_{\Omega }\psi d\gamma _{k}\rightarrow \int_{\Omega }\psi d\gamma $ for
any $\psi \in C_{c}^{2}(\Omega )$, we see from (\ref{gamk}) and the density
of $C_{c}^{2}(\Omega )$ in $C_{0}(\Omega )$ that $(\gamma _{k})$ is weak*
convergent to $\gamma $ on $\Omega $. Further, $\gamma _{k}\geq Nm|_{\Omega
} $ and $\gamma \geq Nm|_{\Omega }$, as in \S \ref{ineq}.

Since $(\gamma -Nm)(\Omega )>0$ by assumption, we can choose a compact set $%
K\subset \Omega $ such that $\alpha >0$, where $\alpha =(\gamma
-Nm)(K^{\circ })$. It follows that $(\gamma _{k}-Nm)(K^{\circ })\geq \alpha
/2$ for all sufficiently large $k$. Let $\gamma _{k}^{\ast }$ denote the
sweeping of $(\gamma _{k}-Nm)|_{K^{\circ }}$ onto $\partial \Omega $. Then
there exists $b>0$ such that $\gamma _{k}^{\ast }\geq b\sigma $ for all
sufficiently large $k$.

If we first consider partial balayage in $O$, then 
\begin{equation}
\mathcal{B}(Nm|_{\Omega }+b\sigma )=Nm|_{U}  \label{bal}
\end{equation}%
for some domain $U$ containing $\overline{\Omega }$. Further, if we choose $b
$ sufficiently small, then $\overline{U}\subset D_{0}$, by Lemma \ref{Lsig}.
By definition 
\begin{equation*}
G_{O}(Nm|_{\Omega }+b\sigma )\geq G_{O}(Nm|_{U})\text{ with equality in }%
O\setminus U.
\end{equation*}%
Since $\overline{U}\subset D_{0}\subset O_{k}\subset O$, this implies that 
\begin{equation*}
G_{O_{k}}(Nm|_{\Omega }+b\sigma )\geq G_{O_{k}}(Nm|_{U})\text{ with equality
in }O_{k}\setminus U.
\end{equation*}%
Thus (\ref{bal}) remains valid if we henceforth consider partial balayage in 
$O_{k}$.

Now let $\eta _{k}$ be the measure $\eta $ constructed as in \S \ref{ineq}.
Since 
\begin{equation*}
\gamma _{k}\geq Nm|_{\Omega }+(\gamma _{k}-Nm)|_{K^{\circ }},
\end{equation*}%
we see from  (\ref{ineqB}) that 
\begin{align*}
\eta_k &=\mathcal{B}((\Delta w_{k,\varepsilon_k})m|_{O_k}) \geq \mathcal{B}((\Delta w_{k,\varepsilon_k})m|_{\Omega}) = \mathcal{B}(\gamma_k)\\&\geq \mathcal{B}(Nm|_{\Omega }+(\gamma
_{k}-Nm)|_{K^{\circ }}) =\mathcal{B}(Nm|_{\Omega }+\gamma _{k}^{\ast })\\ &\geq \mathcal{B}%
(Nm|_{\Omega }+b\sigma )=Nm|_{U},
\end{align*}
By choosing a suitable subsequence of $(\delta _{k})$ we can arrange that $%
(G_{O_{k}}\eta _{k})$ converges locally uniformly on $U$ to a function $u$
in $C^{1}(U)$ satisfying $\Delta u=N$ in $U$, and then $(\nabla
G_{O_{k}}\eta _{k})$ also converges locally uniformly on $U$ to $\nabla u$.
From the final inequality in (\ref{bound}) we see that $\left\Vert \nabla
u\right\Vert \leq \lim_{k\rightarrow \infty }r_{D_{k}}=r_{\Omega }$ on $U$,
so the proposition is proved.

\bigskip

We can now easily address the case of equality in Theorem \ref{main}. We
first assume that $\lambda _{1}(\Omega )=r_{\Omega }$ and choose $v$ as in
Theorem \ref{brenier}, with $D=\Omega $. We claim that $\Delta v\equiv N$ on 
$\Omega $. Indeed, if this were not the case, then there would exist$\ u,U$
as in the above proposition. Since $r_{\Omega }\leq \left\Vert \nabla
u\right\Vert _{\partial \Omega }$, the subharmonic function $\left\Vert
\nabla u\right\Vert $ would then achieve its maximum inside $U$, forcing it
to be constant. Thus $\left\Vert \mathrm{Hess}(u)\right\Vert ^{2}=\Delta
(\left\Vert \nabla u\right\Vert ^{2})=0$, contradicting the fact that $%
\Delta u=N$ on $U$. Hence $\Delta v=N$ on $\Omega $. Since the Hessian of $v$
has determinant equal to $1$, all its eigenvalues are $1$, and so it is the
identity matrix. Thus $\nabla v$ is a translation, and $\Omega $ is a ball
(of radius $r_{\Omega }$).

Conversely, if $\Omega $ is a ball of radius $r$, then (as we noted in the
introduction) $\lambda _{1}(\Omega )=\lambda (\Omega )=r$.

\section{Minimizers for $\protect\lambda _{1}(\Omega )$\label{minimiz}}

In general we do not know whether there exist minimizers $u$ for the
definition of $\lambda _{1}(\Omega )$ in (\ref{L1}). However, we can give
separate necessary and sufficient conditions for a function $u$ in $C^{1}(%
\overline{\Omega })$ to be a minimizer when such minimizers exist. We begin
with a necessary condition.

\begin{proposition}
Suppose that $u\in C^{1}(\overline{\Omega })$ and $\Delta u=N$ in $\Omega $.
If $\Vert \nabla u\Vert _{\overline{\Omega }}=\lambda _{1}(\Omega ),$ then $%
\Vert \nabla u\Vert $ is constant on $\partial \Omega $ (whence $\left\Vert
\nabla u\right\Vert =\lambda _{1}(\Omega )$ on $\partial \Omega $, by (\ref%
{gradsubh})).
\end{proposition}

\bigskip

\noindent \textbf{Proof. }Let $\lambda _{1}=\lambda _{1}(\Omega )$. Suppose
that {$u\in C^{1}(\overline{\Omega })$ satisfies $\Delta u=N$ and }$\Vert
\nabla u\Vert _{\overline{\Omega }}=\lambda _{1}(\Omega )${,} and that there
exists $y_{0}\in \partial \Omega $ such that $\left\Vert \nabla
u(y_{0})\right\Vert <\lambda _{1}$. Then we can choose an open ball $B_{0}$
centred at $y_{0}$ such that 
\begin{equation}
a:=\left\Vert \nabla u\right\Vert _{B_{0}\cap \overline{\Omega }}<\lambda
_{1}.  \label{a}
\end{equation}%
We choose $\Omega _{0}$ to be a domain with $C^{2}$ boundary such that $%
\Omega \subset \Omega _{0}\subset \Omega \cup B_{0}$ and $\Omega
_{0}\backslash \overline{\Omega }\neq \emptyset $.

Let $f=-G_{\Omega }(Nm)$, let $g$ be a continuous extension of $\left. 
\dfrac{\partial f}{\partial n}\right\vert _{\partial \Omega \cap \partial
\Omega _{0}}$ to $\partial \Omega _{0}$ and $0<\varepsilon <\lambda _{1}$.
Proposition 4 of Sakai \cite{Sak} tells us that there is a finite sum $\mu $
of point masses (of variable sign) in $\Omega _{0}\backslash \overline{%
\Omega }$ satisfying%
\begin{equation*}
\left\vert \frac{\partial }{\partial n}G_{\Omega _{0}}\mu +g\right\vert
<\varepsilon \text{ \ \ on \ }\partial \Omega _{0}.
\end{equation*}%
Since 
\begin{eqnarray*}
G_{\Omega _{0}}\mu &\in &C^{1}(\overline{\Omega }_{0}\backslash \mathrm{supp}%
\mu ), \\
\nabla G_{\Omega _{0}}\mu &=&\left( -\frac{\partial }{\partial n}G_{\Omega
_{0}}\mu \right) n\text{ \ \ on \ }\partial \Omega _{0}, \\
\nabla f &=&\frac{\partial f}{\partial n}n\text{ \ \ on \ }\partial \Omega ,
\end{eqnarray*}%
there exists $\kappa >0$ such that 
\begin{equation}
\left\Vert \nabla f-\nabla G_{\Omega _{0}}\mu \right\Vert <\varepsilon \text{
\ on \ }U_{\kappa }\text{,}  \label{gf}
\end{equation}%
where 
\begin{equation*}
U_{\kappa }=\{x\in \Omega :\mathrm{dist}(x,\partial \Omega \cap \partial
\Omega _{0})<\kappa \}.
\end{equation*}

Now let $\delta >0$ and%
\begin{equation*}
u_{\delta }=(1-\delta )u+\delta \left( f-G_{\Omega _{0}}\mu \right) \text{ \
\ on \ }\overline{\Omega }.
\end{equation*}%
Clearly $u_{\delta }\in C^{1}(\overline{\Omega })$, $\Delta u_{\delta }=N$
on $\Omega $ and%
\begin{equation*}
\left\Vert \nabla u_{\delta }\right\Vert \leq (1-\delta )\left\Vert \nabla
u\right\Vert +\delta \left\Vert \nabla f-\nabla G_{\Omega _{0}}\mu
\right\Vert .
\end{equation*}%
Hence, by (\ref{gf}),%
\begin{equation}
\left\Vert \nabla u_{\delta }\right\Vert <(1-\delta )\left\Vert \nabla
u\right\Vert +\delta \varepsilon \leq (1-\delta )\lambda _{1}+\delta
\varepsilon <\lambda _{1}\text{ \ \ on \ }U_{\kappa }.  \label{ua}
\end{equation}%
Also, 
\begin{equation}
\left\Vert \nabla u_{\delta }\right\Vert \leq a+\delta \left\Vert \nabla
f-\nabla G_{\Omega _{0}}\mu \right\Vert _{B_{0}\cap \Omega }\text{ \ \ on \ }%
B_{0}\cap \Omega  \label{ub}
\end{equation}%
and, in view of (\ref{a}), the right hand side above can be made less than $%
\lambda _{1}$ by choosing $\delta $ to be sufficiently small. Combining (\ref%
{ua}) and (\ref{ub}), we obtain a contradiction to the hypothesis that $u$
achieves the minimum value  in (\ref{L1}). Hence $\left\Vert
\nabla u\right\Vert =\lambda _{1}$ on $\partial \Omega $, as claimed.

\bigskip

The converse to the above proposition is false, as we will now explain. Let us recall that $\lambda _{1}(B(0,1))=1$. We claim that, when 
$N\geq 3$, there are functions $u\in C^{1}(\overline{B(0,1)})$\ satisfying\ $%
\Delta u=N$ on $B(0,1)$ and $\left\Vert \nabla u\right\Vert =c$ on $\partial
B(0,1)$, yet $\lambda _{1}(B(0,1))\neq c$. For example, if we define 
\begin{equation*}
u(x)=\frac{N}{2(N-2)}\left( \sum_{i=1}^{N-1}x_{i}^{2}-x_{N}^{2}\right) ,
\end{equation*}%
then $\Delta u=N$ on $B(0,1)$ and $\left\Vert \nabla u\right\Vert =N/(N-2)$
on $\partial B(0,1)$.

Similarly, if $N=2$ and $E=\{(x,y):4x^{2}+y^{2}<1\}$, then the function $%
u(x,y)=2x^{2}-y^{2}$ satisfies $\Delta u=2$ on $E$ and $\left\Vert \nabla
u\right\Vert =2$ on $\partial E$. However, as we will see later, the actual
value of $\lambda _{1}(E)$ is $2/3$.

\bigskip

Next we give a sufficient condition for a function $u\in C^{1}(\overline{%
\Omega })$ to be a minimizer for the definition of $\lambda _{1}(\Omega )$
in (\ref{L1}).

\begin{proposition}
\label{suff}Let $u\in C^{1}(\overline{\Omega })$,\ where\ $\Delta u=N$ on $%
\Omega $. If $\left\Vert \nabla u\right\Vert $ is constant on $\partial
\Omega $ and $\nabla u\cdot n\geq 0$ on $\partial \Omega $, then $\lambda
_{1}(\Omega )=\left\Vert \nabla u\right\Vert _{\overline{\Omega }}$.
Further, any two functions satisfying these hypotheses differ only by a
constant.
\end{proposition}

\noindent \textbf{Proof. }Suppose that $\left\Vert \nabla u\right\Vert =c$
and $\nabla u\cdot n\geq 0$ on $\partial \Omega $. Let $v\in C^{1}(\overline{%
\Omega })$, where $\Delta v=N$ on $\Omega $, and define $w=v-u$. Then $w\in 
\mathcal{H}$. We choose a point $y\in \partial \Omega $ at which $w$
achieves its maximum value. If $w$ is non-constant, then the Hopf boundary
point lemma (see Section 6.4.2 of Evans \cite{Ev}) tells us that $\nabla
w(y)\cdot n_{y}>0$ and so $\nabla w(y)$ is actually a positive multiple of $%
n_{y}$. Since $\nabla u\cdot n\geq 0$ on $\partial \Omega $, we now see that%
\begin{equation*}
\left\Vert \nabla v\right\Vert _{\overline{\Omega }}\geq \left\Vert \nabla
v(y)\right\Vert >\left\Vert \nabla u(y)\right\Vert =c.
\end{equation*}%
It follows that $c=\lambda _{1}(\Omega )$, as required.

Finally, the preceding argument shows that any two functions satisfying the
hypotheses differ only by a constant.

\bigskip

A further useful sufficient condition for a function $u$ to be a minimizer
in (\ref{L1}) applies when $u$ is locally convex.

\begin{theorem}
\label{lconv} Let $u\in C^{1}(\overline{\Omega })$, where $u$ is locally
convex on $\Omega $ and satisfies $\Delta u=N$ there. If $\Vert \nabla
u\Vert =c$ on $\partial \Omega $, then \newline
(i) $\lambda _{1}(\Omega )=c$;\newline
(ii) $\nabla u:\Omega \rightarrow B(0,\lambda _{1}(\Omega ))$ is surjective.
\end{theorem}

\noindent \textbf{Proof. }Let $v=\Vert \nabla u\Vert ^{2}$. Then $\Delta
v=2\left\Vert \mathrm{Hess}(u)\right\Vert ^{2}\geq 0$ on $\Omega $. In fact, 
$\Delta v>0$ on a dense open subset of $\Omega $, for otherwise $\nabla u$
would be constant on a nonempty open set, contradicting the hypothesis that $%
\Delta u=N$. Let $\Omega _{\varepsilon }=\{x\in \Omega :\Vert \nabla
u(x)\Vert <c-\varepsilon \}$, where $\varepsilon >0$ is sufficiently small
that $\Omega _{\varepsilon }\neq \emptyset $. Clearly $\overline{\Omega }%
_{\varepsilon }\subset \Omega $. Further, by the Hopf boundary point lemma, $%
\partial v/\partial n>0$ on $\partial \Omega _{\varepsilon }$, so $n=\nabla
v/\Vert \nabla v\Vert $ on $\partial \Omega _{\varepsilon }$. Since $\mathrm{%
Hess}(u)$ is positive semidefinite on $\Omega $, 
\begin{equation*}
\nabla u\cdot n=\nabla u\cdot \frac{\nabla v}{\Vert \nabla v\Vert }=\frac{2}{%
\Vert \nabla v\Vert }\sum_{i=1}^{N}\sum_{j=1}^{N}\frac{\partial ^{2}u}{%
\partial x_{i}\partial x_{j}}\frac{\partial u}{\partial x_{i}}\frac{\partial
u}{\partial x_{j}}\geq 0\text{ \ \ on }\partial \Omega _{\varepsilon }.
\end{equation*}%
Given any $x\in \partial \Omega $, we choose $B(y,r)\subset \mathbb{R}%
^{N}\backslash \overline{\Omega }$ such that $\partial B(y,r)\cap \overline{%
\Omega }=\{x\}$ and then $x_{\varepsilon }\in \partial \Omega _{\varepsilon
} $ satisfying $\left\Vert x_{\varepsilon }-y\right\Vert =\mathrm{dist}(y,%
\overline{\Omega }_{\varepsilon })$ to see that 
\begin{equation*}
\left( \nabla u\cdot n\right) (x)=\frac{y-x}{\Vert y-x\Vert }\cdot \nabla
u(x)=\lim_{\varepsilon \rightarrow 0}\frac{y-x_{\varepsilon }}{\Vert
y-x_{\varepsilon }\Vert }\cdot \nabla u(x_{\varepsilon })\geq 0,
\end{equation*}%
so part (i) follows from Proposition \ref{suff}. We note, for future
reference, that 
\begin{equation}
\lambda _{1}(\Omega _{\varepsilon })=\lambda _{1}(\Omega )-\varepsilon .
\label{eps}
\end{equation}

We will now show, further, that $\nabla u\cdot n>0$ on $\partial \Omega
_{\varepsilon }$. We write $x=(x^{\prime },x_{N})\in \mathbb{R}^{N-1}\times 
\mathbb{R}$ and choose our co-ordinate system so that $0\in \partial \Omega
_{\varepsilon }$, that the normal $n_{0}$ is in the direction of $(0,...,0,1)
$, and that the Hessian of the function $x^{\prime }\mapsto u(x^{\prime },0)$
at $0^{\prime }$ is a diagonal matrix. We know from the proof of part (i)
that $(\partial u/\partial x_{N})(0)\geq 0$. Now suppose, for the sake of
contradiction, that $(\partial u/\partial x_{N})(0)=0$. Since $v$ is
constant on $\partial \Omega _{\varepsilon }$ we see that 
\begin{equation*}
\frac{\partial v}{\partial x_{i}}(0)=0\text{, \ whence }\ \frac{\partial u}{%
\partial x_{i}}(0)\frac{\partial ^{2}u}{\partial x_{i}^{2}}(0)=0\text{ \ \ \ 
}(i=1,...,N-1).
\end{equation*}%
We reorder the first $N-1$ coordinates so that, for some $m\in \{1,...,N\}$,%
\begin{equation}
\frac{\partial u}{\partial x_{i}}(0)=0\text{ \ \ }(1\leq i\leq m-1)\text{ \
\ and \ \ }\frac{\partial ^{2}u}{\partial x_{i}^{2}}(0)=0\text{ \ \ }(m\leq
i\leq N-1),  \label{zero}
\end{equation}%
and define 
\begin{equation*}
a_{N}=\frac{\partial ^{2}u}{\partial x_{N}^{2}}(0)\text{ \ \ and \ \ }b_{i}=%
\frac{\partial ^{2}u}{\partial x_{i}\partial x_{N}}(0)\text{ \ \ }(1\leq
i\leq N-1).
\end{equation*}%
By (\ref{zero}) the Hessian of the function $(x_{m},...,x_{N})\mapsto
u(0,...,0,x_{m},...,x_{N})$ at $(0,...,0)$ has the form

\begin{equation*}
\left( 
\begin{array}{ccccc}
0 & 0 & \cdots & 0 & b_{m} \\ 
0 & 0 & \cdots & 0 & b_{m+1} \\ 
\vdots & \vdots & \vdots & \vdots & \vdots \\ 
0 & 0 & \cdots & 0 & b_{N-1} \\ 
b_{m} & b_{m+1} & \cdots & b_{N-1} & a_{N}%
\end{array}%
\right) .
\end{equation*}%
Hence $b_{i}=0$ when $m\leq i\leq N-1$, because this submatrix of $\mathrm{%
Hess}(u)$ is also positive semidefinite. By (\ref{zero}) and the Hopf
boundary point lemma, we now arrive at the contradiction%
\begin{equation*}
0<\frac{\partial v}{\partial x_{N}}(0)=2\sum_{i=1}^{N}\frac{\partial u}{%
\partial x_{i}}(0)\frac{\partial ^{2}u}{\partial x_{i}\partial x_{N}}(0)=0.
\end{equation*}%
Thus 
\begin{equation}
\nabla u\cdot n>0\text{ \ on \ }\partial \Omega _{\varepsilon }\text{ \ \ \ }%
(\varepsilon >0).  \label{pos}
\end{equation}

We will now establish (ii). Let $y\in \partial B(0,1)$ and define 
\begin{equation*}
u_{t}(x)=u(x)-ty\cdot x\text{ \ \ \ }(x\in \Omega ;0\leq t<\lambda
_{1}(\Omega _{\varepsilon })).
\end{equation*}%
Then $\nabla u_{t}=\nabla u-ty$ and $\Delta u_{t}=N$ in $\Omega
_{\varepsilon }$. Let 
\begin{equation*}
A=\{t\in \lbrack 0,\lambda _{1}(\Omega _{\varepsilon })):\text{there exists }%
x\in \Omega _{\varepsilon }\text{ such that }\nabla u_{t}(x)=0\}.
\end{equation*}%
It follows from (\ref{pos}) that $u$ cannot attain the value $\min_{%
\overline{\Omega }_{\varepsilon }}u$ on $\partial \Omega _{\varepsilon }$,
so $0\in A$. We will now show that $A$ is both open and closed relative to $%
[0,\lambda _{1}(\Omega _{\varepsilon }))$. To see that $A$ is closed, let $%
(t^{(k)})$ be a sequence in $A$ that converges to some $t\in \lbrack
0,\lambda _{1}(\Omega _{\varepsilon }))$. There exist points $x^{(k)}$ in $%
\Omega _{\varepsilon }$ such that $\nabla u_{t^{(k)}}(x^{(k)})=0$ and, by
choosing a subsequence, we may arrange that $(x^{(k)})$ converges to some
point $x\in \overline{\Omega }_{\varepsilon }$. Clearly $\nabla u_{t}(x)=0$,
so $x\notin \partial \Omega _{\varepsilon }$, because 
\begin{equation*}
\left\Vert \nabla u_{t}\right\Vert \geq \left\Vert \nabla u\right\Vert
-t=\lambda _{1}(\Omega )-\varepsilon -t=\lambda _{1}(\Omega _{\varepsilon
})-t>0\text{ \ \ on \ \ }\partial \Omega _{\varepsilon }
\end{equation*}%
by (\ref{eps}). Hence $A$ is closed. To see that $A$ is open, let $t\in A$,
choose $x\in \Omega $ such that $\nabla u_{t}(x)=0$, and define 
\begin{equation*}
\Omega ^{\prime }=\{z\in \Omega _{\varepsilon }:\Vert \nabla u_{t}(z)\Vert
<\alpha \},\text{ \ \ where \ \ }\alpha =\inf_{\partial \Omega _{\varepsilon
}} \Vert \nabla u_{t}\Vert.
\end{equation*}%
Then $\alpha >0$ and $x\in \Omega ^{\prime }$. We can apply the result of
the previous paragraph to $u_{t}$ to see that $\nabla u_{t}\cdot n>0$ \ on \ 
$\partial \Omega ^{\prime }$. When $\left\vert s\right\vert $ is
sufficiently small, the function $u_{t+s}$ thus also has a strictly positive
normal derivative on $\partial \Omega ^{\prime }$, and so attains the value $%
\min_{\overline{\Omega ^{\prime }}}u_{t+s}$ at some point $x^{(s) }\in
\Omega ^{\prime }$. Since $\nabla u_{t+s}(x^{(s) })=0$, we see that the set $%
A$ is also open relative to $[0,\lambda _{1}(\Omega _{\varepsilon }))$.

Hence $A=[0,\lambda _{1}(\Omega _{\varepsilon }))$. It follows that, for any 
$z\in B(0,\lambda _{1}(\Omega _{\varepsilon }))$, there exists $x\in \Omega
_{\varepsilon }$ such that $\nabla u(x)=z$, and so $\nabla u:\Omega
_{\varepsilon }\rightarrow B(0,\lambda _{1}(\Omega _{\varepsilon }))$ is
surjective. Finally, we let $\varepsilon \rightarrow 0+$ and note from (\ref%
{eps}) that $\lambda _{1}(\Omega _{\varepsilon })\rightarrow \lambda
_{1}(\Omega )$.

\bigskip

\noindent \textbf{Example. }Consider the ellipsoid 
\begin{equation*}
E=\left\{ x\in \mathbb{R}^{N}:\sum_{i=1}^{N}a_{i}^{2}x_{i}{}^{2}<1\right\} ,
\end{equation*}%
where $a_{i}>0$ $(i=1,...,N)$. The function%
\begin{equation*}
u(x)=\frac{N}{2\sum_{i=1}^{N}a_{i}}\sum_{i=1}^{N}a_{i}x_{i}^{2}
\end{equation*}%
is clearly convex, satisfies $\Delta u=N$ on $E$ and $\left\Vert \nabla
u(x)\right\Vert =N/\sum_{i=1}^{N}a_{i}$ on $\partial E$. Thus 
\begin{equation*}
\lambda _{1}(E)=\frac{N}{a_{1}+a_{2}+\dots +a_{N}}.
\end{equation*}

\bigskip

\bigskip

Finally, we remark that minimizers need not be locally convex functions. For
example, if $N\geq 3$ and $\Omega =B(0,R)\backslash \overline{B(0,r)}$,
where $R>r>0$, then the function 
\begin{equation*}
u(x)=\frac{1}{2}\Vert x\Vert ^{2}+\frac{R+r}{(N-2)(R^{1-N}+r^{1-N})}\Vert
x\Vert ^{2-N}\quad (x\in \Omega )
\end{equation*}%
satisfies $\Delta u=N$ on $\Omega $\ and\ $\nabla u=cn$\ on\ $\partial
\Omega $, where%
\begin{equation*}
c=\frac{R^{N}-r^{N}}{R^{N-1}+r^{N-1}}.
\end{equation*}%
Thus $\lambda _{1}(\Omega )=c$, by Proposition \ref{suff}. However, along
the $x_{N}$-axis, $\mathrm{Hess}(u)$ is a diagonal matrix in which the first 
$N-1$ diagonal entries are valued%
\begin{equation*}
1-\frac{R+r}{R^{1-N}+r^{1-N}}x_{N}^{-N},
\end{equation*}%
so $u$ is not locally convex near the inner boundary. (An analogous example
when $N=2$ may be obtained by replacing $\Vert x\Vert ^{2-N}/(N-2)$ with $%
\log (1/\Vert x\Vert )$ in the formula for $u(x)$ above.)

\bigskip

\noindent \textit{Stephen J. Gardiner and Marius Ghergu}

\noindent School of Mathematics and Statistics

\noindent University College Dublin

\noindent Dublin 4, Ireland

\noindent stephen.gardiner@ucd.ie

\noindent marius.ghergu@ucd.ie

\bigskip

\noindent \textit{Tomas Sj\"{o}din}

\noindent Department of Mathematics

\noindent Link\"{o}ping University

\noindent 581 83, Link\"{o}ping

\noindent Sweden

\noindent tomas.sjodin@liu.se

\end{document}